\documentclass[11pt]{article}
\usepackage{amsmath,amssymb,amsfonts,amsthm,graphicx}
\usepackage{hyperref}
\newtheorem{theorem}{Theorem}[section]
\newtheorem{proposition}[theorem]{Proposition}
\newtheorem{corollary}[theorem]{Corollary}
\newtheorem{lemma}[theorem]{Lemma}
\newtheorem{remark}[theorem]{Remark}
\newtheorem{definition}[theorem]{Definition}
\numberwithin{equation}{section}
\title{A GAUSS-KUZMIN THEOREM AND RELATED QUESTIONS FOR $\theta$-EXPANSIONS}
\author{
    Gabriela Ileana Sebe\footnote{e-mail: igsebe@yahoo.com.} \\
    \emph{\small Politehnica University of Bucharest, Faculty of Applied Sciences},\\
    \emph{\small Splaiul Independentei 313, 060042, Bucharest, Romania} \\
    and\\    
    Dan Lascu\footnote{e-mail: lascudan@gmail.com.}\nonumber \\
    \emph{\small Mircea cel Batran Naval Academy, 1 Fulgerului, 900218 Constanta,
    Romania} \\
    }
\begin{document}
\maketitle
\thispagestyle{empty}
\begin{abstract}
Using the natural extension for $\theta$-expansions, we give an infinite-order-chain
representation of the sequence of the incomplete quotients of these expansions. 
Together with the ergodic behavior of a certain homogeneous random system with complete connections, 
this allows us to solve a variant of Gauss-Kuzmin problem for the above fraction expansion.
\end{abstract}
{\bf Mathematics Subject Classifications (2010).}  \\
{\bf Key words}: continued fractions, $\theta$-expansions, Perron-Frobenius operator, random system with complete connections, Gauss-Kuzmin problem. 

\section{Introduction}

During the last fifty years a large amount of research has been devoted to the study of various algorithms for the representation of real numbers by means of sequences of integers. Motivated by problems in random continued fraction expansions (see \cite{BG-2000}), Chakraborty and Rao \cite{CR-2003} have initiated a systematic study of the continued fraction expansion of a number in terms of an irrational $\theta \in (0,1)$. This new expansion of positive reals, different from the regular continued fraction expansion is called $\theta$\textit{-expansion.}

The purpose of this paper is to solve a Gauss-Kuzmin problem for $\theta$-expansions. 
In order to solve the problem, we apply the theory of random systems with complete connections extensively studied by Iosifescu and Grigorescu \cite{IG-2009}.
First we outline the historical framework of this problem. 
In Section 1.2, we present the current framework. 
In Section 1.3, we review known results. 

\subsection{Gauss' Problem}

One of the first and still one of the most important results in the metrical theory of continued fractions is the so-called Gauss-Kuzmin theorem.
Any irrational $0<x<1$ can be written as the infinite regular continued fraction 

\begin{equation}
x = \displaystyle \frac{1}{a_1+\displaystyle \frac{1}{a_2+\displaystyle \frac{1}{a_3+ \ddots}}} :=[a_1, a_2, a_3, \ldots], \label{1.1}
\end{equation}
where $a_n \in \mathbb{N}_+ : = \left\{1, 2, 3, \ldots\right\}$ \cite{IK-2002}.
Such integers $a_1, a_2, \ldots$ are called {\it incomplete quotients} (or {\it continued fraction digits}) of $x$.
The metrical theory of continued fraction expansions started on 25th October 1800, with a note by Gauss in his mathematical diary \cite{Brez}. 
Define the \textit{regular continued fraction} (or \textit{Gauss}) \textit{transformation} $\tau$ on the unit interval $I:=[0, 1]$ by 
\begin{equation}
\tau (x) = \left\{\begin{array}{lll} 
\displaystyle \frac{1}{x}-\left\lfloor \displaystyle \frac{1}{x} \right\rfloor & \hbox{if} & x \neq 0, \\ 
\\
0 & \hbox{if} & x = 0,
\end{array} \right. \label{1.2}
\end{equation}
where $\left\lfloor \cdot \right\rfloor$ denotes the floor (or entire) function.
With respect to the asymptotic behavior of iterations $\tau^{n}=\tau\circ \cdots\circ \tau$ ($n$-times) of $\tau$, Gauss wrote (in modern notation) that
\begin{equation}
\lim_{n \rightarrow \infty} \lambda \left(\tau^n \leq x\right) = \frac{\log (1+x)}{\log 2},
\quad x \in I, \label{1.3}
\end{equation}
where $\lambda$ denotes the Lebesgue measure on $I$.
In 1812, Gauss asked Laplace \cite{Brez} to estimate the {\it $n$-th error term} $e_n(x)$ defined by
\begin{equation}
e_n(x) := \lambda (\tau^{-n}[0, x]) - \frac{\log (1+x)}{\log 2}, \quad n \geq 1, \ x\in I. \label{1.4}
\end{equation}
This has been called \textit{Gauss' Problem}. 
It received first solution more than a century later, when R.O. Kuzmin \cite{Kuzmin-1928} showed in 1928 that $e_n(x) = \mathcal{O}(q^{\sqrt{n}})$ as $n \rightarrow \infty$, uniformly in $x$ with some (unspecified) $0 < q < 1$. 
This has been called the {\it Gauss-Kuzmin theorem} or the {\it Kuzmin theorem}.

One year later, using a different method, Paul L\'evy \cite{Levy-1929} improved Kuzmin's result by showing that $\left|e_n(x)\right| \leq q^n$ for $n \in \mathbb{N}_+$, $x \in I$, with $q = 3.5 - 2\sqrt{2} = 0.67157...$. 
For such historical reasons, the {\it Gauss-Kuzmin-L\'evy theorem} 
is regarded as the first basic result in the rich metrical theory of continued fractions. 
An advantage of  the Gauss-Kuzmin-L\'evy theorem
relative to the Gauss-Kuzmin theorem is the determination of the value of $q$.

To this day the Gauss transformation, on which metrical theory of regular continued fraction is based, has fascinated researchers from various branches of mathematics and science with many applications in computer science, cosmology and chaos theory \cite{C-1992}. In the last century, mathematicians broke new ground in this area. Apart from the regular continued fraction expansion, very many other continued fraction expansions were studied \cite{RS-1992,Schweiger}. 

By such a development,
generalizations of these problems 
for non-regular continued fractions 
are also called as
the {\it Gauss-Kuzmin problem}
and the {\it Gauss-Kuzmin-L\'evy problem} \cite{IS-2006,L-2013,Sebe-2000,Sebe-2001,Sebe-2002}.

\subsection{$\theta$-expansions}

For a fixed $\theta \in (0,1)$, we start with a brief review of continued fraction expansion with respect to $\theta$, analogous to the regular continued fraction expansion which corresponds to the case $\theta=1$.

For $x \in (0,\infty)$ let 
\[
a_0 := \max\{n \geq 0: n\theta \leq x \}. 
\]
If $x$ equals $a_0 \theta$, we write 
\[
x := [a_0 \theta]. 
\]
If not, define $r_1$ by 
\[
x := a_0 \theta + \frac{1}{r_1} 
\]
where $0 < 1/r_1 < \theta$. Then 
$
r_1 > 1/\theta \geq \theta 
$
and let 
\[
a_1 := \max \{n \geq 0: n \theta \leq r_1\}. 
\]
If $r_1 = a_1 \theta$, then we write 
\[
x := [a_0 \theta, a_1 \theta], 
\]
i.e., 
\[
x = a_0 \theta + \frac{1}{a_1 \theta}. 
\]
If $a_1 \theta < r_1$, define $r_2$ by 
\[
r_1 := a_1 \theta + \frac{1}{r_2} 
\]
where $0 < 1/r_2 < \theta$. So, $r_2 > 1/\theta \geq \theta$ and let 
\[
a_2 := \max \{n \geq 0: n\theta \leq r_2\}. 
\]
In this way, either the process terminates after a finite number of steps or it continues indefinitely.
Following standard notation, in the first case we write 
\begin{equation}
x := [a_0 \theta; a_1 \theta,\ldots, a_n \theta] \label{1.5'}
\end{equation}
and we call this \textit{the finite continued fraction expansion of} $x$ \textit{with respect to} $\theta$ (\textit{terminating at the} $n$\textit{-stage}). In the second case, we write 
\begin{equation}
x := [a_0 \theta; a_1 \theta, a_2 \theta, \ldots] \label{1.5''}
\end{equation}
and we call this \textit{the infinite continued fraction expansion of} $x$ \textit{with respect to }$\theta$. 

When $0 < x < \theta$, we have $a_0 = 0$ and instead of writing 
\begin{equation}
x := [0; a_1 \theta, a_2 \theta, \ldots], \label{1.5'''}
\end{equation}
we simply write 
\begin{equation}
x := [a_1 \theta, a_2 \theta, \ldots] \label{1.5''''}
\end{equation}
which is the same in the usual notation
\begin{equation}
x = \frac{1}{\displaystyle a_1\theta
+\frac{1}{\displaystyle a_2\theta 
+ \frac{1}{\displaystyle a_3\theta + \ddots} }}. \label{1.5}
\end{equation}
Such $a_n$'s are also called \textit{incomplete quotients} (or \textit{continued fraction digits}) of $x$ with respect to the expansion in $(\ref{1.5})$.

In general, the $\theta$-expansion of a number $x > 0$ is
\begin{equation}
a_0 \theta + [0; a_1 \theta, a_2 \theta, \ldots] := [a_0 \theta; a_1 \theta, a_2 \theta, \ldots] \label{1.6}
\end{equation}
where $a_0 = \left\lfloor x/ \theta\right\rfloor$.

For $x \in [0, \theta]$, the $\theta$-continued fraction expansion of $x$ in (\ref{1.5}) leads to an analogous transformation of Gauss map $\tau$ in (\ref{1.2}). A natural question is whether this new transformation admits an absolutely continuous invariant probability like the Gauss measure in the case $\theta = 1$. Until now, the invariant measure was identified only in the particular case $\theta^2 = 1/m$, $m$ a positive integer \cite{CR-2003}.

Motivated by this argument, since the invariant measure is a crucial tool in our approach, in the sequel we will consider only the case $\theta^2 = 1/m$ with $m$ a positive integer. Then $[a_1\theta, a_2\theta, a_3\theta, \ldots]$ is the $\theta$-expansion of any $x \in [0, \theta]$ if and only if the following conditions hold: 
\begin{enumerate}
\item [(i)]
$a_n \geq m$ for any $m \in \mathbb{N}_+$;

\item [(ii)]
in case when $x$ has a finite expansion, i.e., $x = [a_1\theta, a_2\theta, a_3\theta, \ldots, a_n\theta]$, then $a_n \geq m+1$.
\end{enumerate}

This continued fraction is treated as the following dynamical system.

\begin{definition} \label{def.1.1}
Let $\theta \in (0,1)$ and $m \in \mathbb{N}_+$ such that $\theta^2 = 1/m$. 
\begin{enumerate}

\item [(i)]
The measure-theoretical dynamical system $([0,\theta],{\cal B}_{[0,\theta]}, T_{\theta})$ is defined as follows:
$\mathcal{B}_{[0,\theta]}$ denotes the $\sigma$-algebra of all Borel subsets of $[0,\theta]$, and $T_{\theta}$ is the transformation 
\begin{equation}
T_{\theta}: [0,\theta] \to [0,\theta];\quad 
T_{\theta}(x):=
\left\{
\begin{array}{ll}
{\displaystyle \frac{1}{x} - \theta \left \lfloor \frac{1}{x \theta} \right\rfloor} & 
{\displaystyle \mbox{if } x \in (0, \theta],}\\
\\
0 & \mbox{if } x=0.
\end{array}
\right. \label{1.7}
\end{equation} 

\item [(ii)] 
In addition to (i), we write $([0,\theta], {\cal B}_{[0,\theta]}, \gamma_{\theta}, T_{\theta} )$ as $( [0,\theta],{\cal B}_{[0,\theta]},T_{\theta})$ with the following probability measure $\gamma_{\theta}$ on $( [0,\theta], {\cal B}_{[0,\theta]})$:

\begin{equation}
\gamma_{\theta} (A) := 
\frac{1}{\log \left(1+\theta^{2}\right)}
\int_{A} \frac{\theta dx}{1 + \theta x}, 
\quad A \in {\mathcal{B}}_{[0,\theta]}. \label{1.8}
\end{equation}

\end{enumerate}
\end{definition}

By using $T_{\theta}$, the sequence $(a_n)_{n \in {\mathbb{N}}_+}$ in (\ref{1.5}) is obtained as follows: 
\begin{equation}
a_n=a_n(x) = a_1\left(T_{\theta}^{n-1}(x)\right), \quad n \in {\mathbb{N}}_+, \label{1.9}
\end{equation}
with $T_{\theta}^0 (x) = x$ and
\begin{equation}
a_1=a_1(x) = \left\{\begin{array}{lll} 
\lfloor \frac{1}{x \theta} \rfloor & \hbox{if} & x \neq 0, \\ 
\infty & \hbox{if} & x = 0.
\end{array} \right. \label{1.10}
\end{equation}

In this way, $T_{\theta}$ gives the algorithm of $\theta$-expansion.

\begin{proposition} \label{prop.1.2}
Let 
$([0,\theta],{\cal B}_{[0,\theta]},\gamma_{\theta},T_{\theta})$ 
be as in Definition \ref{def.1.1}(ii).
\begin{enumerate}
\item [(i)]
$([0,\theta],{\cal B}_{[0,\theta]},\gamma_{\theta},T_{\theta})$ is ergodic.
\item [(ii)]
The measure $\gamma_{\theta}$ is invariant under $T_{\theta}$,
that is, 
$\gamma_{\theta} (A) = \gamma_{\theta} (T^{-1}_{\theta}(A))$ 
for any $A \in {\mathcal{B}}_{[0, \theta]}$.
\end{enumerate}
\end{proposition}
\noindent \textbf{Proof.} 
See Section 8 in \cite{CR-2003}.
\hfill $\Box$\\

By Proposition \ref{prop.1.2}(ii), $([0,\theta],{\cal B}_{[0,\theta]},\gamma_{\theta},T_{\theta})$ is a 
``dynamical system" in the sense of Definition 3.1.3 in \cite{BG}.

\subsection{Known results and applications}
In this subsection we recall known results and their applications for $\theta$-expansions. 

\subsubsection{Known results} 

Let $0 < \theta < 1$. Define the {\it $n$-th order convergent} $[a_1\theta, a_2 \theta, \ldots, a_n \theta]$ of $x \in [0, \theta]$ by truncating the $\theta$-expansion in (\ref{1.5}). Thus, Chakraborty and Rao proved in \cite{CR-2003} that
\begin{equation}
[a_1 \theta, a_2 \theta, \ldots, a_n \theta] \to x, \quad n\to \infty. \label{1.12}
\end{equation}

In what follows the stated identities hold for all $n$ in case $x$ has an infinite $\theta$-expansion and they hold for $n \leq k$ in case $x$ has a finite $\theta$-expansion terminating at the $k$-th stage.

To this end, define real functions $p_n(x)$ and $q_n(x)$, for $n \in \mathbb{N}_+$, by 
\begin{eqnarray}
p_n(x) &:=& a_n(x) \theta p_{n-1}(x) + p_{n-2}(x), \quad \label{1.13} \\
q_n(x) &:=& a_n(x) \theta q_{n-1}(x) + q_{n-2}(x), \quad \label{1.14}
\end{eqnarray}
with $p_{-1}(x) := 1$, $p_0(x) := 0$, $q_{-1}(x) := 0$ and $q_{0}(x) := 1$. By induction, we have
\begin{equation}
p_{n-1}(x)q_{n}(x) - p_{n}(x)q_{n-1}(x) = (-1)^{n}, \quad n \in \mathbb{N} := \mathbb{N}_+ \cup \{0\}. \label{1.15}
\end{equation}
By using (\ref{1.13}) and (\ref{1.14}), we can verify that 
\begin{equation}
x = \frac{p_n(x) + T^n_{\theta}(x)p_{n-1}(x)}
{q_n(x) + T^n_{\theta}(x)q_{n-1}(x)}, \quad n \geq 1. \label{1.16}
\end{equation}
By taking $T^n_{\theta}(x)=0$ in (\ref{1.16}), we obtain $ [a_1\theta, a_2\theta, \ldots, a_n\theta] = p_n(x)/q_n(x)$. Using (\ref{1.15}) and (\ref{1.16}) we obtain 
\begin{equation}
x - \frac{p_n(x)}{q_n(x)} = \frac{(-1)^{n+1} T^n_{\theta}(x)}{q_n(x)(q_n(x)+T^n_{\theta}(x)q_{n-1}(x))}, \quad n \geq 1. \label{1.17}
\end{equation}
By applying $0 \leq T^n_{\theta} \leq \theta$ to (\ref{1.17}), we can verify that
\begin{equation}
\frac{1}{q_n(x)(q_{n+1}(x)+\theta q_n(x))} \leq \left|x - \frac{p_n(x)}{q_n(x)}\right| \leq \frac{1}{q_n(x)q_{n+1}(x)}. \label{1.18}
\end{equation}
From (\ref{1.14}), we have that $q_n(x) \geq \theta$, $n \in \mathbb{N}_+$. Further, also from (\ref{1.14}) and by induction we have that
\begin{equation}
q_n(x) \geq \left\lfloor \frac{n}{2}\right\rfloor \theta^2. \label{1.19}
\end{equation}
Finally, (\ref{1.12}) follows from (\ref{1.18}) and (\ref{1.19}).

\subsubsection{Application to ergodic theory}

Similarly to classical results on regular continued fractions, using the ergodicity of $T_{\theta}$ and Birkhoff's ergodic theorem \cite{DK-2002}, a number of results were obtained.

For $q_n$ in (\ref{1.14}), its asymptotic growth rate $\beta$ is defined as 
\begin{equation}
\beta = \lim_{n \rightarrow \infty} \frac{1}{n} \log q_n. \label{1.20}
\end{equation}
This is a L\'evy result and Chakraborty and Rao \cite{CR-2003} obtained that $\beta$ is a finite number 
\begin{equation}
\beta = \frac{1}{1 + \theta^2} \int_{0}^{\theta} \frac{\theta \log x}{1 + x \theta}\mathrm{d}x. \label{beta}
\end{equation}
They also give a Khintchin result, i.e., the asymptotic value of the arithmetic mean of $a_1, a_2, \ldots, a_n$ where $a_1$ and $a_n$ are given in (\ref{1.10}) and (\ref{1.9}). We have
\begin{equation}
\lim_{n \rightarrow \infty} \frac{a_1 + a_2 + \ldots + a_n}{n} = \infty. \label{1.21}
\end{equation}

It should be stressed that the ergodic theorem does not yield any information on the convergence rate in the Gauss problem that amounts to the asymptotic behavior of $\mu(T_{\theta}^{-n})$ as $n \rightarrow \infty $, where $\mu$ is an arbitrary probability measure on ${\cal B}_{[0,\theta]}$. 

It is only very recently that there has been any investigation of the metrical properties of the $\theta$-expansions. Thus,  the results obtained in this paper allow to a solution of a Gauss-Kuzmin type problem. We may emphasize that, to our knowledge, Theorem \ref{G-K-L} is the first Gauss-Kuzmin result proved for $\theta$-expansions. Our solution presented here is based on the ergodic behavior of a certain random system with complete connections.

The paper is organized as follows.
In Section 2, we show the probability structure of $(a_n)_{n \in {\mathbb{N}}_+}$ 
under the Lebesgue measure by using the Brod\'en-Borel-L\'evy formula. 
In Section 3, we consider the so-called natural extension 
of $([0, \theta],{\cal B}_{[0, \theta]},\gamma_{\theta},T_{[0, \theta]})$ \cite{Nakada}.
In Section 4, we derive its Perron-Frobenius operator under different probability measures on $([0, \theta],{\cal B}_{[0, \theta]})$. 
Especially, 
we derive the asymptotic behavior
for the Perron-Frobenius operator of $(([0, \theta],{\cal B}_{[0, \theta]},\gamma_{\theta},T_{[0, \theta]})$. 
In Section 5, we study the ergodicity of the associated random system with complete connections (RSCC for short).
In Section 6, we solve a variant of Gauss-Kuzmin problem for $\theta$-expansions. 
By using the ergodic behavior of the RSCC introduced in Section 5,
we determine the limit of the sequence 
$(\,\mu (T_{\theta}^n < x)\,)_{n\geq 1}$ of distributions
as $n \rightarrow \infty$.

\section{Prerequisites}

Roughly speaking, the metrical theory of continued fraction expansions is about the sequence 
$(a_n)_{n \in {\mathbb{N}}_+}$ of incomplete quotients and related sequences \cite{IK-2002}.
As remarked earlier in the introduction we will adopt a similar strategy to that used for regular continued fractions. We begin with a Brod\'en-Borel-L\'evy formula for $\theta$-expansions. Then some consequences of it to be used in the sequel are also derived.


In this section let us fix $0 < \theta < 1$, $\theta^2 = 1/m$, $m \in \mathbb{N}_+$. 

For $x \in [0, \theta]$ consider $a_n=a_n(x)$, $n \in \mathbb{N}_+$, as in (\ref{1.9}) and (\ref{1.10}). Putting $\mathbb{N}_m := \{m, m+1,\ldots\}$, $m \in \mathbb{N}_+$, the incomplete quotients $a_n$, $n \in \mathbb{N}_+$, take positive integer values in $\mathbb{N}_m$. 

For any $n \in \mathbb{N}_+$ and $i^{(n)}=(i_1, \ldots, i_n) \in \mathbb{N}_m^n$, define the {\it fundamental interval associated with} $i^{(n)}$ by
\begin{equation}
I (i^{(n)}) = \{x \in [0, \theta]:  a_k(x) = i_k \mbox{ for } k=1, \ldots, n \}, \label{2.01}
\end{equation}
where $I (i^{(0)}) = [0, \theta]$. 
For example, for any $i \in \mathbb{N}_m$ we have
\begin{equation}
I\left(i\right) = \left\{x \in [0, \theta]: a_1(x) = i \right\} = \left( \frac{1}{(i+1)\theta}, \frac{1}{i \theta} \right). \label{2.02}
\end{equation}

We will write $I(a_1, \ldots, a_n) = I\left(a^{(n)}\right)$, $n \in \mathbb{N}_+$. If $n \geq 1$ and $i_n \in \mathbb{N}_m$, then we have
\begin{equation}
I(a_1, \ldots, a_n) = I\left(i^{(n)}\right). \label{2.03}
\end{equation}

From the definition of $T_{\theta}$ and (\ref{1.16}), we have 
\begin{equation}
I(a^{(n)}) = (u(a^{(n)}), v(a^{(n)})), \label{2.04}
\end{equation}
where $u(a^{(n)})$ and $v(a^{(n)})$ are defined as
\begin{equation}
u\left(a^{(n)}\right) := \left\{
\begin{array}{lll}
	\displaystyle\frac{p_n+ \theta p_{n-1}}{q_n+ \theta q_{n-1}} & \quad \mbox{if $n$ is odd}, \\
	\\
	\displaystyle\frac{p_n}{q_n} & \quad \mbox{if $n$ is even}, \\ 
\end{array}
\right. \label{2.05}
\end{equation}
and
\begin{equation}
v\left(a^{(n)}\right) := \left\{
\begin{array}{lll}
	\displaystyle\frac{p_n}{q_n} & \quad \mbox{if $n$ is odd}, \\
	\\
	\displaystyle\frac{p_n+ \theta p_{n-1}}{q_n+\theta q_{n-1}} & \quad \mbox{if $n$ is even},
\end{array}
\right. \label{2.06}
\end{equation}
where $p_n:=p_n(x)$ and $q_n:=q_n(x)$ are defined in (\ref{1.13}) and (\ref{1.14}), respectively.

Let $\lambda_{\theta}$ denote a Lebesgue measure on $[0, \theta]$. Using (\ref{1.15}) we get
\begin{eqnarray}
\lambda_{\theta}\left(I\left(a^{(n)}\right)\right) &=& \frac{1}{\theta} \left| \frac{p_n}{q_n} - \frac{p_n+ \theta p_{n-1}}{q_n+\theta q_{n-1}} \right| \nonumber \\
&=& \frac{1}{q_n (q_n + \theta q_{n-1})}. \label{2.07}
\end{eqnarray}

To derive the so-called Brod\'en-Borel-L\'evy formula \cite{IG-2009, IK-2002} for $\theta$-expansions, let us define $(s_n)_{n \in {\mathbb{N}}}$ by
\begin{equation}
s_0 := 0,\quad 
s_n := q_{n-1}/q_{n}, \, \quad n \geq 1. \label{2.14}
\end{equation}
From (\ref{1.14}), $s_n = 1/(a_n \theta + s_{n-1})$ for $n \geq 1$. Hence 
\begin{equation}
s_n = \frac{1}{\displaystyle a_n\theta
+\frac{1}{\displaystyle a_{n-1}\theta 
+ \ddots + \frac{1}{a_1 \theta}}} := [a_n\theta, a_{n-1}\theta, \ldots, a_1\theta], \label{2.15}
\end{equation} 
for $n \geq 1$. 

\begin{proposition} [Brod\'en-Borel-L\'evy-type formula] \label{prop.BBL}
Let $\lambda_{\theta}$ denote the Lebesgue measure on $[0, \theta]$.
For any $n \in \mathbb{N}_+$, 
the conditional probability 
$\lambda_{\theta} (T^n_{\theta} < x |a_1,\ldots, a_n )$ is given as follows:
\begin{equation}
\lambda_{\theta} (T^n_{\theta} < x |a_1,\ldots, a_n ) 
= \frac{(s_n \theta + 1)x}{\theta(s_n x+1)}, \quad x \in [0, \theta], \label{2.16}
\end{equation}
where $s_n$ is defined in $(\ref{2.14})$ and $a_1, \ldots, a_n$ are as in $(\ref{1.9})$ and $(\ref{1.10})$.
\end{proposition}

\noindent\textbf{Proof.} 
By definition, we have
\begin{equation}
\lambda_{\theta} (T^n_{\theta} < x |a_1,\ldots, a_n ) = \frac{\lambda_{\theta}\left(\left(T^n_{\theta} < x\right) \cap I(a_1,\ldots, a_n) \right)}{\lambda_{\theta}\left(I(a_1,\ldots, a_n)\right)} \label{2.17}
\end{equation}
for any $n \in \mathbb{N}_+$ and $x \in [0, \theta]$.
Using (\ref{1.16}) and (\ref{2.04}) we get
\begin{eqnarray}
\lambda_{\theta}\left(\left(T^n_{\theta} < x \right) \cap I(a_1,\ldots, a_n)\right) &=& \frac{1}{\theta}\left|\frac{p_n}{q_n} - \frac{p_n+xp_{n-1}}{q_n+xq_{n-1}}\right| \nonumber \\
\nonumber \\
&=& \frac{x}{ q_n (q_n + x q_{n-1})\theta}. \nonumber
\end{eqnarray}
From this and (\ref{2.07}) it follows that
\begin{eqnarray}
\lambda_{\theta}\left(T^n_{\theta} < x |a_1, \ldots, a_n \right) &=& \frac{\lambda_{\theta}\left(\left(T^n_{\theta} < x\right) \cap I(a_1,\ldots, a_n) \right)}{\lambda_{\theta}\left(I(a_1,\ldots, a_n)\right)} \nonumber \\  
\nonumber \\
&=& \frac{x\left(q_n + \theta q_{n-1}\right)}{(q_n +x q_{n-1})\theta}= \frac{x(s_n \theta + 1)}{(s_n x+1)\theta}, \label{2.18}
\end{eqnarray}
for any $n \in \mathbb{N}_+$ and $x \in [0, \theta]$. 
\hfill $\Box$

The Brod\'en-Borel-L\'evy formula allows us to determine the probability structure of incomplete quotients $(a_n)_{n \in {\mathbb{N}}_+}$ under $\lambda_{\theta}$.

\begin{proposition} \label{prop.2.3}
For any $i \in \mathbb{N}_m$ and $n \in \mathbb{N}_+$, we have
\begin{equation}
\lambda_{\theta}(a_1 = i) = \frac{m}{i(i+1)}, \quad 
\lambda_{\theta}\left(a_{n+1}=i |a_1,\ldots, a_n \right) = P_i(s_n), \label{2.19}
\end{equation}
where $(s_n)_{n \in {\mathbb{N}}_+}$ is defined in (\ref{2.14}), and 
\begin{equation}
P_i(x) := \frac{x \theta + 1}{(x + i\theta)(x + (i+1)\theta )}. \label{2.20}
\end{equation}
\end{proposition}
\noindent\textbf{Proof.}
From (\ref{2.02}), the case $\lambda_{\theta}(a_1=i)$ holds. For $n \in \mathbb{N}_+$ and $x \in [0, \theta]$,
we have $ T_{\theta}^n(x) = [a_{n+1} \theta, a_{n+2} \theta, \ldots]$.
Using (\ref{2.16}) we obtain
\begin{eqnarray}
\lambda_{\theta}(\,a_{n+1}=i \,|\,a_1,\ldots, a_n\, ) 
& = & \lambda_{\theta}\left(\,T^n_{\theta} \in \left(\frac{1}{(i+1)\theta}, \frac{1}{i\theta}\right]\, 
| \,a_1,\ldots, a_n \,\right). \nonumber \\
&=& \frac{(s_n \theta + 1)\frac{1}{i\theta}}{\theta(s_n \frac{1}{i\theta}+1)} 
- \frac{(s_n \theta + 1)\frac{1}{(i+1)\theta}}{\theta(s_n \frac{1}{(i+1)\theta}+1)} \nonumber \\
\nonumber\\
&=& P_i(s_n).  \label{2.21}
\end{eqnarray}
\hfill $\Box$\\
\begin{remark}

\begin{enumerate}

\item[(i)]
It is easy to check that 
\begin{equation}
\sum_{i=m}^{\infty}P_i(x)=1\quad \mbox{ for any } x \in [0, \theta]. \label{2.9}
\end{equation}
\item[(ii)]
Proposition \ref{prop.2.3} is the starting point of an approach to the metrical theory of $\theta$-expansions via dependence with complete connections (see \cite{IG-2009}, Section 5.2)
\end{enumerate}

\end{remark}
\begin{corollary}
The sequence $(s_n)_{n \in {\mathbb{N}}_+}$ with $s_0=0$ is an $[0, \theta]$-Markov chain on $([0,\theta],{\cal B}_{[0,\theta]}, \lambda_{\theta})$ with the following transition mechanism: from state $s$ the possible transitions are to any state $1/(s + i \theta)$ with corresponding transition probability $P_i(s)$, $i \in \mathbb{N}_m$.
\end{corollary}

\section{An infinite-order-chain representation}

In this section we introduce the natural extension $\overline{T_{\theta}}$ of $T_{\theta}$ in (\ref{1.7}) and define 
extended random variables according to Chap.1.3.3 of \cite{IK-2002}. Then we give an infinite-order-chain representation of the sequence of the incomplete quotients for $\theta$-expansions. 

\subsection{Natural extension}

Let $([0,\theta],{\cal B}_{[0,\theta]}, T_{\theta})$ be as in Definition 1.1(i). 
Let be $[0, \theta]^{2} := [0, \theta] \times [0, \theta]$ and the square space $([0, \theta]^2,{\cal B}^2_{[0, \theta]}) := ([0, \theta],{\cal B}_{[0, \theta]})\times ([0, \theta],{\cal B}_{[0, \theta]})$. 
\begin{definition}
The natural extension \cite{Nakada}
of $([0, \theta],{\cal B}_{[0, \theta]},T_{\theta})$
is
$([0, \theta]^2,{\cal B}^2_{[0, \theta]},\overline{T_{\theta}})$ 
where the transformation $\overline{T_{\theta}}$ 
of $([0, \theta]^2,{\cal B}^2_{[0, \theta]})$ is defined by 
\begin{equation}
\overline{T_{\theta}}:[0, \theta]^2\to [0, \theta]^2; \quad
\overline{T_{\theta}}(x,y) := \left(T_{\theta}(x), \, \frac{1}{a_1(x)\theta + y}\right), \quad (x, y) \in [0, \theta]^2. \label{3.3}
\end{equation}
\end{definition}

\noindent
This is a one-to-one transformation of $[0, \theta]^2$ with the inverse
\begin{equation}
(\overline{T_{\theta}})^{-1}(x, y) 
= 
\left(\frac{1}{a_1(y)\theta + x}, \,
T_{\theta}(y)\right), \quad \,(x, y) \in [0, \theta]^{2}\,. \label{3.4}
\end{equation}
Iterations of (\ref{3.3}) and (\ref{3.4}) are given as follows for each $n \geq 2$:
\begin{equation}
\left(\overline{T_{\theta}}\right)^n(x, y) = 
(\,T_{\theta}^n(x), \,[a_n(x) \theta, a_{n-1}(x) \theta, \ldots, 
a_2(x) \theta,\, a_1(x) \theta + y ] \,), \label{3.5} \\
\end{equation}
\begin{equation}
\left(\overline{T_{\theta}}\right)^{-n}(x, y) = 
(\,[a_n(y)\theta, a_{n-1}(y)\theta, \ldots, a_2(y)\theta, 
\,a_1(y)\theta + x ],\, T^{n}_{\theta}(y) \,). 
\label{3.6}
\end{equation}
For $\gamma_{\theta}$ in (\ref{1.8}), define its {\it extended measure} $\overline{\gamma_{\theta}}$
on $([0, \theta]^2,{\mathcal{B}}^2_{[0, \theta]})$ as
\begin{equation}
\overline{\gamma_{\theta}} (B) := \frac{1}{\log(1 + \theta^2)} \int\!\!\!\int_{B} 
\frac{\mathrm{d}x\mathrm{d}y}{(1 + xy)^2},
\quad B \in {\mathcal{B}}^2_{[0, \theta]}. \label{3.7}
\end{equation}
Then
\begin{equation}
\overline{\gamma_{\theta}}(A \times [0, \theta]) 
= \overline{\gamma_{\theta}}([0, \theta] \times A) = \gamma_{\theta}(A) \label{3.7'}
\end{equation}
for any $A \in {\mathcal{B}}_{[0, \theta]}$. 
\begin{proposition} \label{prop.3.2}
The measure $\overline{\gamma_{\theta}}$ is preserved by $\overline{T_{\theta}}$.
\end{proposition}
\noindent\textbf{Proof.} 
We show that $\overline{\gamma_{\theta}} ((\overline{T_{\theta}})^{-1}(B)) 
= \overline{\gamma_{\theta}} (B)$ for any $B \in {\mathcal{B}}^2_{[0, \theta]}$. 
Since $\overline{T_{\theta}}$ is invertible 
on $[0, \theta]^2$, the last equation is equivalent to 
\begin{equation}
\overline{\gamma_{\theta}} (\overline{T_{\theta}}(B)) = \overline{\gamma_{\theta}} (B)\quad 
 \mbox{for any } B \in {\mathcal{B}}^2_{[0, \theta]}. \label{3.51}
\end{equation}
Recall fundamental interval in (\ref{2.01}).
Since the collection of Cartesian products $I\left(i^{(n)}\right) \times I\left(j^{(t)}\right)$, $i^{(n)} \in \mathbb{N}^n_m$, $j^{(t)} \in
\mathbb{N}^t_m$, $n, t \in \mathbb{N}$, generates the $\sigma$-algebra ${\mathcal{B}}^2_{[0, \theta]}$, it is enough to show that
\begin{equation}
\overline{\gamma_{\theta}} \left(\overline{T_{\theta}}\left(I\left(i^{(n)}\right) \times I\left(j^{(t)}\right)\right)\right) = 
\overline{\gamma_{\theta}} \left(I\left(i^{(n)}\right) \times I\left(j^{(t)}\right)\right) \label{3.8'} 
\end{equation}
for any $i^{(n)} \in \mathbb{N}^n_m$, $j^{(t)} \in \mathbb{N}^t_m$, $n, t \in \mathbb{N}$.
It follows from (\ref{3.7'}) and Proposition \ref{prop.1.2}(ii) that (\ref{3.8'}) holds for $n = 0$ and $t \in \mathbb{N}$. If $n \in \mathbb{N}_+$, then it is easy to see that
\begin{equation}
\overline{T_{\theta}}\left(I\left(i^{(n)}\right) \times I\left(j^{(t)}\right)\right) = 
I(i_2, \ldots ,i_n) \times I(i_1, j_1, \ldots ,j_t), \quad t \in \mathbb{N}_+,
\end{equation}
where $I(i_2, \ldots, i_n)$ equals $[0, \theta]$ for $n = 1$. Also, if $I\left(i^{(n)}\right) = (a, b) \subset [0, \theta]$ and $I\left(j^{(t)}\right) = (c,d) \subset [0, \theta]$, then
\[
I(i_2, \ldots, i_n) = \left(\frac{1}{b} - i_1 \theta, \frac{1}{a} - i_1 \theta\right)
\]
and
\[
I(i_1, j_1, \ldots, j_t) = \left(\frac{1}{d+i_1 \theta}, \frac{1}{c+i_1\theta}\right). \nonumber
\]
A simple computation yields
\begin{equation}
\overline{\gamma_{\theta}} \left((a, b) \times (c,d)\right) = \frac{1}{\log(1+\theta^2)} \log \frac{(ac+1)(bd+1)}{(ad+1)(bc+1)} 
\end{equation}
and then 
\begin{eqnarray}
\overline{\gamma_{\theta}} \left(\left(\frac{1}{b} - i_1 \theta, \frac{1}{a} - i_1 \theta\right) \times \left(\frac{1}{d+i_1 \theta}, \frac{1}{c+i_1\theta}\right)\right) \qquad \qquad \qquad \qquad \qquad \nonumber \\
= \frac{1}{\log(1+\theta^2)} \log \frac{((1/b- i_1 \theta)/(d+i_1 \theta)+1)((1/a -i_1 \theta)/(c+i_1\theta)+1)}
{((1/b-i_1 \theta)/(c+i_1\theta)+1)((1/a-i_1\theta)/(d+i_1\theta)+1)} \nonumber \\
= \frac{1}{\log(1+\theta^2)} \log \frac{(ac+1)(bd+1)}{(ad+1)(bc+1)}, \qquad \qquad \qquad \qquad \qquad \qquad \qquad \qquad 
\end{eqnarray}
that is, (\ref{3.8'}) holds.
\hfill $\Box$

\subsection{Extended random variables}

With respect to $\overline{T_{\theta}}$ in (\ref{3.3}), 
define {\it extended incomplete quotients} $\overline{a}_l(x,y)$, 
$l \in \mathbb{Z}:=\{\ldots, -2, -1, 0, 1, 2, \ldots\}$, $(x, y) \in [0, \theta]^2$ by
\begin{equation}
\overline{a}_{l+1}(x, y) 
:= \overline{a}_1((\overline{T_{\theta}})^{l} (x, y) \,),
\quad l \in \mathbb{Z}, \label{3.8}
\end{equation}
with 
\begin{equation}
\overline{a}_1 (x, y) = a_1(x), \quad (x, y) \in [0, \theta]^2.
\end{equation}
\begin{remark} 

\begin{enumerate}

\item[(i)]
Since $\overline{T_{\theta}}$ is invertible it follows that 
$\overline{a}_{l}(x, y)$ in (\ref{3.8}) 
is also well-defined for $l \leq 0$.
By (\ref{3.5}) and (\ref{3.6}), we have
\begin{equation}
\overline{a}_n(x, y) = a_n(x), \quad 
\overline{a}_0(x, y) = a_1(y), \quad 
\overline{a}_{-n}(x, y) = a_{n+1}(y) \label{3.9}
\end{equation}
for any $n \in \mathbb{N}_+$ and $(x, y) \in [0, \theta]^2$.
\item[(ii)]
From Proposition \ref{prop.3.2}, the doubly infinite sequence 
$(\overline{a}_l(x,y))_{l \in \mathbb{Z}}$ 
is strictly stationary (i.e., its distribution is invariant under a shift of the indices) under $\overline{\gamma}_{\theta}$.
\end{enumerate}

\end{remark}
The following theorem will play a key role in the sequel.

\begin{theorem} \label{th.3.4}
For any $x \in [0, \theta]$ we have 
\begin{equation}
\overline{\gamma_{\theta}} ( [0, x] \times [0, \theta] \,| 
\,\overline{a}_0, \overline{a}_{-1}, \ldots ) 
=\frac{(a \theta +1)x}{(ax + 1)\theta} \quad \overline{\gamma_{\theta}} \mbox{-}\mathrm{a.s.}, \label{3.10}
\end{equation}
where $a:= [\overline{a}_0 \theta, \overline{a}_{-1} \theta, \ldots]$ 
(= the $\theta$-expansion with incomplete quotients $\overline{a}_0$, $\overline{a}_1$, $\ldots$). 
\end{theorem}
\noindent \textbf{Proof.} Let $I_{n}$ denote the fundamental interval $I(\overline{a}_0, \overline{a}_{-1}, \ldots, \overline{a}_{-n})$ for $n \in \mathbb{N}$. We have
\begin{equation}
\overline{\gamma_{\theta}} ( [0, x] \times [0, \theta] \left. \right| \overline{a}_0, \overline{a}_{-1}, \ldots ) = \lim_{n \rightarrow \infty} \overline{\gamma_{\theta}} ( [0, x] \times [0, \theta] \left. \right| \overline{a}_0, \ldots, \overline{a}_{-n} ) \quad \overline{\gamma}_{\theta} \mbox{-a.s.} \label{3.11}
\end{equation}
and
\begin{eqnarray}
\overline{\gamma_{\theta}} ( [0, x] \times [0, \theta] \left. \right| \overline{a}_0, \ldots, \overline{a}_{-n} ) 
&=& 
\displaystyle{\frac{\overline{\gamma_{\theta}} ([0, x] \times I_{n})}{\overline{\gamma_{\theta}} ([0, \theta] \times I_{n})}} \nonumber \\
\nonumber \\
 &=& \displaystyle{\frac{1}{\gamma_{\theta}(I_{n})}\frac{1}{\log(1+\theta^2)}\int_{I_{n}}\mathrm{d}y\displaystyle\int^x_0{\frac{\mathrm{d}u}{(uy+1)^2}}} \nonumber \\
 \nonumber \\
 &=& \displaystyle{\frac{1}{\gamma_{\theta}(I_{n})} \frac{1}{\theta} \int_{I_{n}} \frac{x(1+y \theta)}{1+xy}\, \mathrm{d}\gamma_{\theta}(y)} \nonumber \\
 \nonumber \\
 &=& \frac{x(y_n \theta +1)}{(xy_n+1)\theta}, \label{3.12}
\end{eqnarray}
for some $y_n \in I_{n}$. Since 
\begin{equation}
\lim_{n \rightarrow \infty} y_n =[\overline{a}_0 \theta, \overline{a}_{-1} \theta, \ldots] = a, \label{3.13}
\end{equation}
the proof is completed.
\hfill $\Box$

The probability structure of $(\overline{a}_l)_{l \in \mathbb{Z}}$ under $\overline{\gamma_{\theta}}$ is given as follows. 
\begin{corollary} \label{cor.3.5}
For any $i \in \mathbb{N}_m$, we have
\begin{equation}
\overline{\gamma_{\theta}} (\left.\overline{a}_1 = i\right| \overline{a}_0, \overline{a}_{-1}, \ldots) = P_{i}(a) \quad \overline{\gamma_{\theta}} \mbox{-}\mathrm{a.s.}, \label{3.14}
\end{equation}
where $a = [\overline{a}_0 \theta, \overline{a}_{-1} \theta, \ldots]$ and the functions $P_{i}$, $i \in \mathbb{N}_m$, are defined by $(\ref{2.20})$.
\end{corollary}
\noindent \textbf{Proof.} Let $I_{n}$ be as in the proof of Theorem \ref{th.3.4}. 
We have 
\begin{equation}
\overline{\gamma_{\theta}} (\left.\overline{a}_1 = i\,\right| \,
\overline{a}_0, \overline{a}_{-1}, \ldots) = \lim_{n \rightarrow \infty} 
\overline{\gamma_{\theta}} (\left.\overline{a}_1 = i\,\right| \,I_{n}). \label{3.15}
\end{equation}
We have
\begin{equation}
(\overline{a}_1 = i) = \left(\frac{1}{(i+1)\theta}, \frac{1}{i \theta}\right] \times [0, \theta], \quad i \in \mathbb{N}_m. \label{3.16}
\end{equation}
Now 
\begin{eqnarray}
\overline{\gamma_{\theta}} \left( \left. \left(\frac{1}{(i+1)\theta}, \frac{1}{i \theta}\right] \times [0, \theta] \right| I_{n}\right) &=& 
\frac{\overline{\gamma_{\theta}} \left(\left(\frac{1}{(i+1)\theta}, \frac{1}{i \theta}\right] \times I_{n}\right) }
{\overline{\gamma_{\theta}} ([0, \theta] \times I_{n})} \nonumber \\
\nonumber \\
& = & \frac{1}{\gamma_{\theta} (I_{n})} \int_{I_{n}} P_{i}(y)\, \mathrm{d}\gamma_{\theta}(y) \nonumber \\
\nonumber \\
& = & P_{i}(y_n), \label{3.17}
\end{eqnarray}
for some $y_n \in I_{n}$. From (\ref{3.13}), the proof is completed. 
\hfill $\Box$
\begin{remark} \label{rem.3.6}
The strict stationarity of $\left(\overline{a}_l
\right)_{l \in \mathbb{Z}}$, under $\overline{\gamma_{\theta}}$ implies that
\begin{equation}
\overline{\gamma_{\theta}} 
(\left.\overline{a}_{l+1} = i\, \right|\, \overline{a}_l, 
\overline{a}_{l-1}, \ldots) 
= P_{i}(a) \quad \overline{\gamma}_{\theta} 
\mbox{-}\mathrm{a.s.} \label{3.18}
\end{equation}
for any $i \in \mathbb{N}_m$ and $l \in \mathbb{Z}$, where $a = [\overline{a}_l \theta, \overline{a}_{l-1} \theta, \ldots]$. The last equation emphasizes that $\left(\overline{a}_l\right)_{l \in \mathbb{Z}}$ is an infinite-order-chain in the theory of dependence with complete connections (see \cite{IG-2009}, Section 5.5).
\end{remark}
Define extended random variables $\left(\overline{s}_l\right)_{l \in \mathbb{Z}}$ as
$\overline{s}_l := [\overline{a}_l \theta, \overline{a}_{l-1} \theta, \ldots]$, $l \in \mathbb{Z}$.
Clearly, $\overline{s}_l = \overline{s}_{0} \circ (\overline{T_{\theta}})^l$, $l \in \mathbb{Z}$. It follows from Proposition \ref{prop.3.2} and Corollary \ref{cor.3.5} that $\left(\overline{s}_l\right)_{l \in \mathbb{Z}}$ is a strictly stationary $[0,\theta]$-valued Markov process on 
$([0, \theta]^2,{\mathcal{B}}^2_{[0, \theta]}, \overline{\gamma_{\theta}})$
with the following transition mechanism: from state $\overline{s} \in [0,\theta]$ the possible transitions are to any state $1/(\overline{s} + i \theta)$ with corresponding transition probability $P_i(\overline{s})$, $i \in \mathbb{N}_m$. 
Clearly, for any $l \in \mathbb{Z}$ we have
\begin{equation}
\overline{\gamma_{\theta}}(\overline{s}_l < x ) = \overline{\gamma_{\theta}}( \overline{s}_0 < x) 
 = \overline{\gamma_{\theta}}([0,\theta] \times [0,x)) = \gamma_{\theta}([0,x)), \quad x \in [0, \theta].
\end{equation}

Motivated by Theorem \ref{th.3.4}, we shall consider 
the one-parameter family $\{\gamma_{\theta, a}: a \in [0, \theta]\}$ 
of (conditional) probability 
measures  on $([0, \theta],{\mathcal{B}}_{[0, \theta]})$ 
defined by their distribution functions
\begin{equation}
\gamma_{\theta,a} ([0, x]) := \frac{(a \theta +1)x}{(ax + 1)\theta},
\quad
x, a \in [0, \theta]. \label{3.19}
\end{equation}
Note that $\gamma_{\theta,0} = \lambda_{\theta}$. 

For any $a \in [0, \theta]$ put
\begin{equation}
s_{0,a} := a,\quad
s_{n,a} := \frac{1}{a_n \theta + s_{n-1,a}}, \quad n \in \mathbb{N}_+. \label{3.20}
\end{equation}

\begin{remark}
It follows from the properties just described of the process $(\overline{s}_l)_{l \in \mathbb{Z}}$ that the sequence $(s_{n,a})_{n \in {\mathbb{N}}_+}$ is an $[0, \theta]$-valued Markov chain on 
$([0, \theta],{\cal B}_{[0,\theta]}, \gamma_{\theta,a})$
which starts at $s_{0,a} := a$ and has the following transition mechanism: from state $s \in [0, \theta]$ the possible transitions are to any state $1/(s+i\theta)$ with corresponding transition probability $P_i(s)$, $i \in \mathbb{N}_m$. 
\end{remark}

\section{An operatorial treatment}

Let $([0, \theta],{\cal B}_{[0, \theta]},\gamma_{\theta}, T_{\theta})$ be as in Definition 1.1(ii). 
The Gauss-Kuzmin problem for the transformation $T_{\theta}$ can be approached in terms of the associated Perron-Frobenius operator.

Let $\mu$ be a probability measure on $([0, \theta], {\mathcal{B}}_{[0, \theta]})$ 
such that $\mu((T_{\theta})^{-1}(A)) = 0$ whenever $\mu(A) = 0$ for any
$A \in {\mathcal{B}}_{[0, \theta]}$. 
For example, this condition is satisfied if $T_{\theta}$ is $\mu$-preserving, that is, $\mu (T_{\theta})^{-1} = \mu$. 
Let 
\[
L^1_{\mu}:=\{f: [0, \theta] \rightarrow \mathbb{C} : \int^{\theta}_0 |f |\mathrm{d}\mu < \infty \}.
\]
The {\it Perron-Frobenius operator} 
of \ $T_{\theta}$ under $\mu$ 
is defined as the bounded linear operator $U_{\mu}$  which takes the Banach space $L^1_{\mu}$ into itself and satisfies the equation
\begin{equation}
\int_{A}U_{\mu}f \,\mathrm{d}\mu = \int_{(T_{\theta})^{-1}(A)}f\, \mathrm{d}\mu \quad 
\mbox{ for all }
A \in {\mathcal{B}}_{[0, \theta]},\, f\in L^1_{\mu}. \label{4.1}
\end{equation}
%
%
\begin{proposition} \label{prop.4.1.}
\begin{enumerate}
\item[(i)] 
The Perron-Frobenius operator $U:=U_{\gamma_{\theta}}$ of \ $T_{\theta}$ under the invariant probability measure $\gamma_{\theta}$ is given a.e. in $[0, \theta]$ by the equation
\begin{equation}
Uf(x) = \sum_{i \geq m}P_{i}(x)\,f(u_{i}(x)), \quad m \in \mathbb{N}_+, \, 
f \in L^1_{\gamma_{\theta}}, \label{4.2}
\end{equation}
where $P_{i}$, $i \geq m$, is as in (\ref{2.20}) and $u_i$, $i \geq m$, is defined by
\begin{equation}
u_i : [0, \theta] \rightarrow [0, \theta]; \quad u_i(x) := \frac{1}{x + i \theta}. \label{4.3}
\end{equation}
\item[(ii)]
Let $\mu$ be a probability measure on $([0, \theta],{\mathcal{B}}_{[0, \theta]})$
such that $\mu$ is absolutely continuous with respect to 
the Lebesgue measure $\lambda_{\theta}$ 
and let $h := d\mu / d \lambda_{\theta}$ a.e. in $[0, \theta]$. 
For any $n \in \mathbb{N}$ and $A \in {\mathcal{B}}_{[0, \theta]}$,
we have 
\begin{equation}
\mu \left((T_{\theta})^{-n}(A)\right) 
= \int_{A} U^nf(x) \mathrm{d} \gamma_{\theta}(x), \label{4.4}
\end{equation}
where $f(x):= (\log(1+\theta^2)) \frac{x\theta+1}{\theta}h(x)$, $x \in [0, \theta]$. 
\end{enumerate}

\end{proposition}

\noindent \textbf{Proof.} 
(i) Let $T_{\theta,i}$ denote 
the restriction of $T_{\theta}$ to the subinterval 
$I(i) := \left(\frac{1}{(i+1)\theta}, \frac{1}{i \theta}\right]$, $i \geq m$, that is, 
\begin{equation}
T_{\theta,i}(x) = \frac{1}{x} - i \theta, \quad x \in I(i). \label{4.5}
\end{equation}
Let $C(A):=(T_{\theta})^{-1}(A)$ and $C_{i}(A):=(T_{\theta,i})^{-1}(A)$ 
for any $A \in {\mathcal B}_{[0, \theta]}$. 
Since $C(A)=\bigcup_{i}C_i(A)$ and $C_i\cap C_j$ is a null set when $i \ne j$,
we have
\begin{equation}
\int_{C(A)} f \,\mathrm{d} \gamma_{\theta} 
= \sum_{i \geq m} \int_{C_i(A)}f\, \mathrm{d} \gamma_{\theta},
\quad
f \in L^1_{\gamma_{\theta}},\,A \in {\mathcal{B}}_{[0, \theta]}.
 \label{4.6}
\end{equation}
For any $i \geq m$, by the change of variable
$x = (T_{\theta,i})^{-1}(y) =u_{i}(y)$,
we successively obtain
\begin{eqnarray}
\int_{C_i(A)}f(x) \,\mathrm{d}\gamma_{\theta}(x) &=& \frac{\theta}{\log(1+\theta^2)} \int_{C_i(A)} \frac{f(x)}{1+x\theta}\,\mathrm{d}x \nonumber \\
\nonumber \\
&=& \frac{1}{\log(1+\theta^2)} \int_{A} \frac{f\left(u_{i}(y)\right)}{1+u_{i}(y) \theta} \frac{\theta \mathrm{d}y}{(y+i \theta)^2} \nonumber \\
\nonumber \\
&=& \int_{A} P_{i}(y)\, f\left(u_{i}(y)\right)\,\mathrm{d}\gamma_{\theta}(y). 
\label{4.7}
\end{eqnarray}
Now, (\ref{4.2}) follows from (\ref{4.6}) and (\ref{4.7}).

\noindent
(ii) We will use mathematical induction. For $n=0$, 
the equation (\ref{4.5}) holds by definitions
of $f$ and $h$.
Assume that (\ref{4.5}) 
holds for some $n \in \mathbb{N}$. Then 
\begin{eqnarray}
\mu ((T_{\theta})^{-(n+1)}(A)) &=& \mu ((T_{\theta})^{-n}((T_{\theta})^{-1}(A))) \nonumber \\
\nonumber \\
&=& \int_{C(A)} U^n f(x) \,\mathrm{d} \gamma_{\theta}(x). \label{4.9}
\end{eqnarray}
By the very definition of the Perron-Frobenius operator $U$ we have
\begin{equation}
\int_{C(A)} U^n f(x) \,\mathrm{d}\gamma_{\theta}(x) = \int_{A} U^{n+1} f(x) \,\mathrm{d}\gamma_{\theta}(x). \label{4.10}
\end{equation}
Therefore, 
\begin{equation}
\mu \left((T_{\theta})^{-(n+1)}(A)\right) = \int_{A} U^{n+1} f(x) \mathrm{d}\gamma_{\theta}(x) \label{4.11}
\end{equation}
which ends the proof.

\hfill $\Box$

Let $B([0, \theta])$ denote the collection of all bounded measurable functions $f:[0, \theta] \rightarrow \mathbb{C}$. 
A different interpretation is available for the operator $U$ restricted to $B([0, \theta])$, which is a Banach space under the supremum norm.

\begin{proposition} \label{prop.4.2} 
The operator $U : B([0, \theta]) \rightarrow B([0, \theta])$ is the transition operator of both the Markov chain $(s_{n,a})_{n \in \mathbb{N}_+}$ on $([0, \theta], {\mathcal{B}}_{[0, \theta]}, \gamma_{\theta,a})$, for any $a \in [0, \theta]$, and the Markov chain $(\overline{s}_{l})_{l \in \mathbb{Z}}$ on $([0, \theta]^2, {\mathcal{B}}^2_{[0, \theta]}, \overline{\gamma_{\theta}})$.
\end{proposition} 
\noindent \textbf{Proof.} 
The transition operator of $(s_{n,a})_{n \in \mathbb{N}_+}$ takes $f \in B([0, \theta])$ to the function defined by 
\begin{equation}
E_a\left( \left. f(s_{n+1,a})\right| s_{n,a} = s \right) = \sum_{i \geq m}P_{i}(s)f(u_{i}(s)) = U f(s) \quad \mbox{for any } s \in [0, \theta], \label{4.12}
\end{equation}
where $E_a$ stands for the mean-value operator with respect to the probability measure $\gamma_{\theta,a}$, whatever $a \in [0, \theta]$. 

\hfill $\Box$

A similar reasoning is valid for the case of the Markov	chain $(\overline{s}_l)_{l \in \mathbb{Z}}$.
\begin{remark} \label{rem.4.3}
In hypothesis of Proposition \ref{prop.4.1.}(ii) it follows that 
\begin{equation}
\mu(T_{\theta}^{-n}(A)) - \gamma_{\theta}(A) = \int_{A}(U^{n}f(x)-1)\mathrm{d}\gamma_{\theta}(x),
\end{equation}
for any $n \in \mathbb{N}$ and $A \in {\mathcal{B}}_{[0, \theta]}$, where 
$f(x):= (\log(1+\theta^2))\frac{x\theta+1}{\theta}h(x)$, $x \in [0, \theta]$.
The last equation shows that the asymptotic behavior of $\mu(T_{\theta}^{-n}(A)) - \gamma_{\theta}(A)$ as $n \rightarrow \infty$ is given by the asymptotic behavior of the $n$-th power of the Perron-Frobenius $U$ on $L^1_{\gamma_{\theta}}$ or on smaller Banach spaces.
\end{remark}
\section{Ergodicity of the associated RSCC}

The facts presented in the previous sections lead us to a certain random system with complete connections associated with the $\theta$-expansion. To study the ergodicity of this RSCC it becomes necessary to recall some definitions and results from \cite{IG-2009}.

According to the general theory we have the following statement.
\begin{definition} \label{5.1}
An homogeneous RSCC is a quadruple 
$\left\{(W, {\mathcal W}), (X, {\mathcal X}), u, P\right\}$
where
\begin{enumerate}
	\item[(i)] $(W, {\mathcal W})$ and $(X, {\mathcal X})$ are arbitrary measurable spaces;
	\item[(ii)] $u: W \times X \rightarrow W$ is a $({\mathcal W} \otimes {\mathcal X}, {\mathcal W})$-measurable map;
	\item[(iii)] $P$ is a transition probability function from $(W, {\mathcal W})$ to $(X, {\mathcal X})$.
\end{enumerate}
\end{definition}
For any $n \in \mathbb{N}_+$, consider the maps $u^{(n)} : W \times X^{n} \rightarrow W$, defined by
\[
\left.
\begin{array}{lll} 
u^{(1)}\left(w, x\right) &:=& u(w, x), \\
u^{(n+1)}\left(w, x^{(n+1)}\right) &:=& 
u\left(u^{(n)}\left(w, x^{(n)}\right), x_{n+1}\right), \,  n \geq 1, 
\end{array} \right. \label{5.2}
\]
where $x^{(n)} = (x_1, \ldots, x_n) \in X^{n}$. We will simply write $wx^{(n)}$ for $u^{(n)}(w,x^{(n)})$.
For every $w \in W$, $r \in \mathbb{N}_+$ and $A \in {\mathcal X}^r$, define 
\begin{equation}
\left.
\begin{array}{lll}
P_1(w, A):=P(w, A), \\ 
P_r(w, A):=\displaystyle \int_{X}P(w, \mathrm{d}x_1) \displaystyle\int_{X}P(wx_1, \mathrm{d}x_2)\ldots \displaystyle\int_{X}P(wx^{(r-1)}, \mathrm{d}x_r) \chi_A(x^{(r)}), \, r \geq 2,
\end{array} \right. \label{5.1'}
\end{equation}
where $\chi_A$ is the indicator function of the set $A$.
Obviously, for $n \in \mathbb{N}_+$ fixed, $P_r$ is a transition probability function from $(W, {\mathcal W})$ to $(X^r, {\mathcal X}^r)$.

By virtue of the existence theorem (\cite{IG-2009}, Theorem 1.1.2), 
for a given RSCC $\left\{(W, {\mathcal W}), (X, {\mathcal X}), u, P\right\}$ 
there exists an associated Markov chain with the transition operator $U$ defined by
\[
Uf(w) := \int_{X} P(w, \mathrm{d}x)f(wx), \quad f \in B(W, {\mathcal W}),
\]
where $B(W, {\mathcal W})$ is the Banach space of all bounded ${\mathcal W}$-measurable complex-valued functions defined on $W$.
Moreover, the transition probability function of the associated Markov chain is
\[
Q(w, B) := \int_{X} P(w, \mathrm{d}x) \chi_B(wx) = P(w, B_w), 
\]
where $B_w = \left\{ x \in X: wx \in B \right\}$, $w \in W$, $B \in {\mathcal W}$.
The iterates of the operator $U$ are given by
\[
U^nf(w) = \int_{X^n} P_n(w, \mathrm{d}x^{(n)})f(wx^{(n)}), \quad f \in B(W, {\mathcal W}), n \in \mathbb{N}_+. 
\]
It follows that the $n$-step transition probability function is given by
\[
Q^n(w, B) = P_n(w, B^{(n)}_w), \quad w \in W, \, B \in {\mathcal W}, n \in \mathbb{N}_+,
\]
where $B^{(n)}_w = \left\{ x^{(n)}: wx^{(n)} \in B \right\}$.
Hence the transition operator associated with the Markov chain with state space $(W, {\mathcal W})$ and transition probability function $Q$ is defined by 
\begin{equation}
Uf(\cdot) := \int_{W} Q(\cdot, \mathrm{d}w)f(w), \quad f \in B(W, {\mathcal W}).
\end{equation}
Its iterates are given by
\begin{equation}
U^{n}f (\cdot) = \int_{W} Q^{n}(\cdot, \mathrm{d}w)f(w), \quad n \in \mathbb{N}_+,
\end{equation}
where $Q^n$ is the $n$-step transition probability function. 

Putting 
\[
Q_n(w, B) = \frac{1}{n} \sum_{k=1}^{n} Q^k(w, B)
\]
for all $n \in \mathbb{N}_+$, $w \in W$ and $B \in {\mathcal W}$, it is clear that $Q_n$ is a transition probability function on $(W, {\mathcal W})$. Let $U_n$ be the Markov operator associated with $Q_n$.
Let $(W, d)$ be a metric space and let $L(W)$ denote the Banach space of all complex-valued Lipschitz continuous functions on $W$ with the following norm:
\begin{equation}
\left\| f \right\|_L := \|f\| + s(f), \label{5.3}
\end{equation}
where 
\begin{equation}
\| f \|:= \sup_{w \in W} \left|f(w)\right|, \, 
s(f) := \sup_{w'\neq w''} \frac{|f(w') - f(w'')|}{d(w', w'')}. \label{5.4}
\end{equation}
\begin{definition}
\begin{enumerate}
\item[(i)] 
The operator $U$ is said to be orderly with respect to $L(W)$ if and only if there exists a bounded linear operator $U^\infty$ on $L(W)$ such that
\[
\lim_{n \rightarrow \infty} \left\|U_n -U^\infty \right\|_L = 0. 
\]
\item[(ii)] 
The operator $U$ is said to be aperiodic with respect to $L(W)$ if and only if there exists a bounded linear operator $U^\infty$ on $L(W)$ such that
\[
\lim_{n \rightarrow \infty} \left\|U^n - U^\infty \right\|_L = 0.
\]
\item[(iii)] 
The operator $U$ is said to be ergodic with respect to $L(W)$ if and only if  it is orderly and the range $U^\infty(L(W))$ is one-dimensional.
\item[(iv)] 
The operator $U$ regular with respect to $L(W)$ if and only if it is ergodic and aperiodic.
\end{enumerate}
\end{definition}
\begin{definition}
The transition operator $U$ of a Markov chain with state space $W$ is said to be a Doeblin-Fortet operator
if and only if $U$ takes $L(W)$ into $L(W)$ boundedly with respect to $\|\cdot\|_L$ and there exist $k \in \mathbb{N}_+$, $r \in [0, 1)$ and $R < \infty$ such that 
\begin{equation}
s(U^{k}f) \leq r s(f) + R \| f \|, \quad f \in L(W). \label{5.5}
\end{equation}
Alternatively, the Markov chain itself is said to be a \textit{Doeblin-Fortet chain}.
\end{definition}
The definition below isolates a class of RSCCs, called \textit{RSCCs with contraction}, 
for which the associated Markov chains are Doeblin-Fortet chains. 
\begin{definition} \label{def.contraction}
An RSCC $\left\{\left(W, {\mathcal{W}}\right), \left(X, {\mathcal{X}}\right), u, P\right\}$ 
is said to be with contraction if and only if 
$(W, d)$ is a separable metric space, $r_1 < \infty$, $R_1 < \infty$ and there exists $j \in \mathbb{N}_+$
such that $r_j < 1$. Here 
\begin{equation}
r_j := \sup_{w' \neq w''} \int_{X^j} P_j(w', \mathrm{d}x^{(j)}) \frac{d(w'x^{(j)}, w''x^{(j)})}{d(w', w'')} 
\end{equation}
and 
\begin{equation}
R_j := \sup_{A \in {\mathcal{X}}^{j}} s(P_j (\cdot, A)), \quad j \in \mathbb{N}_+
\end{equation}
where $P_j$ is a transition probability function from $\left(W, {\mathcal{W}}\right)$ to 
$\left(X^j, {\mathcal{X}}^{j}\right)$ defined in (\ref{5.1'}).
\end{definition}
We shall also need the following results. 
\begin{theorem} \label{th.DF}
The Markov chain associated with an RSCC with contraction is a Doeblin-Fortet chain. 
\end{theorem}
\begin{lemma}
Assume that the Markov operator $U$ is aperiodic with respect to $L(W)$. 
Put $T = U-U^{\infty}$. Then we have $U^{\infty}U^{\infty}=U^{\infty}$, $TU^{\infty}=U^{\infty}T=0$, $U^n=U^{\infty}+T^n$, $n \in \mathbb{N}_+$. 
Moreover, there exist positive constants $q<1$ and $K$ such that
\begin{equation}
\left\|T^{n}\right\|_L \leq Kq^n, \quad n \in \mathbb{N}_+.
\end{equation}
\end{lemma}
This immediately implies the validity of
\begin{equation}
\left\|U^nf-U^{\infty}f\right\|_L \leq K q^n \left\|f\right\|_L, \quad f \in L(W), n \in \mathbb{N}_+. \label{5.13}
\end{equation}
\begin{definition}
A Markov chain is said to be compact if and only if its state space is a compact metric space 
$(W, d)$ and its transition operator is a Doeblin-Fortet operator.
\end{definition}
\begin{theorem} \label{th.ordered}
A compact Markov chain is orderly with respect to $L(W)$ and there exists a transition probability function $Q^{\infty}$ on $(W,{\mathcal{W}})$ such that
\begin{equation}
   U^{\infty}f(\cdot)= \int_{W} Q^{\infty}(\cdot,\mathrm{d}w) f(w),\quad f \in L(W).
\end{equation}
Moreover, $Q^{\infty}(\cdot,B) \in E(1)$ for any $B \in {\mathcal{W}}$ and for any $w \in W$, $Q^{\infty}(w, \cdot)$ is a stationary probability for the chain. 
(Here $E(1)$ is the set of the eigenvalues of modulus $1$ of the operator $U$).
\end{theorem}
A criterion of regularity for a compact Markov chain is expressed in Theorem \ref{th.5.9}
in terms of the supports $\sigma_n(w)$ of the transition probability functions $Q^{n}(w, \cdot)$, $n \in \mathbb{N}_+$.
\begin{theorem} \label{th.5.9}
A compact Markov chain is regular with respect to $L(W)$ if and only if there exists a point 
$w^{*} \in W$ such that 
\begin{equation}
\lim_{n \rightarrow \infty} d\left(\sigma_n(w), w^{*}\right) = 0, \quad w \in W.
\end{equation}
\end{theorem}
The application of this criterion is facilitated by the inter-relationship among the sets $\sigma_n(w)$, $n \in \mathbb{N}_+$, which is given in the next lemma.
\begin{lemma} \label{lem.supp}
For all $m, n \in \mathbb{N}$ and $w \in W$, we have
\begin{equation}
\sigma_{m+n}(w) = \overline{\bigcup_{w' \in \sigma_m(w)}\sigma_n(w')},
\end{equation}
where the overline mean topological closure in $W$.
\end{lemma}
Now, we are able to study the following RSCC
\begin{equation}
\left\{\left([0, \theta], {\mathcal{B}}_{[0, \theta]}\right), \left(\mathbb{N}_m, {\mathcal{P}}(\mathbb{N}_m)\right), u, P\right\}, \label{5.17}
\end{equation}
where 
$u : [0, \theta] \times \mathbb{N}_m \rightarrow [0, \theta]$, 
$u(s,i)=u_i(s)$ is given in (\ref{4.3}) 
and the function $P(s,i)=P_i(s)$ given in (\ref{2.20}) defines a transition probability from 
$([0,\theta], {\mathcal{B}}_{[0, \theta]})$ to $(\mathbb{N}_m, \mathcal{P}\left((\mathbb{N}_m\right)))$. 
Here $\mathbb{N}_m = \{m, m+1, \ldots\}$, $m \in \mathbb{N}_+$
and $\mathcal{P} \left(\mathbb{N}_m\right)$ denotes the power set of $\mathbb{N}_m$.

Whatever $a \in [0,\theta]$ the Markov chain $(s_{n,a})_{n \in \mathbb{N}}$ associated with the RSCC (\ref{5.1}) has the transition operator $U$, with the transition probability function
\begin{equation}
Q(s,B) = \sum_{ \{ \left.i \geq m \right| u_i(s) \in B \} } P_i(s), \quad s \in [0,\theta], B \in {\mathcal{B}}_{[0, \theta]}.
\end{equation}
Then $Q^{n}(\cdot, \cdot)$ will denote the $n$-step transition probability function of the same Markov chain.

\begin{proposition}
RSCC (\ref{5.17}) is regular with respect to $L([0,\theta])$. 
Moreover there exist a stationary probability measure $Q^{\infty} = \gamma_{\theta}$ and two positive constants $q<1$ and $K$ such that
\begin{equation}
\left\|U^n f- \int_{0}^{\theta} f \mathrm{d}\gamma_{\theta} \right\|_L \leq K q^n \left\|f\right\|_L, \quad n \in \mathbb{N}_+, f \in L([0,\theta]),
\end{equation}
where
\begin{equation}
U^{n}f (\cdot):= \int_{0}^{\theta} Q^{n}(\cdot, \mathrm{d}s)f(s) \label{5.20}
\end{equation}
\end{proposition}
\noindent \textbf{Proof.} 
Since 
\begin{eqnarray}
\frac{\mathrm{d}}{\mathrm{d}s} u_i(s) &=& \frac{-1}{(s+i \theta)^2}, \nonumber \\
\frac{\mathrm{d}}{\mathrm{d}s} P_i(s) &=& \frac{i \theta - 1/\theta}{(s+i \theta)^2} - \frac{(i+1) \theta - 1/\theta}{(s+(i+1) \theta)^2}, \nonumber 
\end{eqnarray}
for any $s \in [0, \theta]$ and $i \geq m$ it follows that 
\[
\sup_{s \in [0, \theta]} \left| \frac{\mathrm{d}}{\mathrm{d}s} u_i(s) \right| = \frac{1}{(i \theta)^2}
\]
and
\[
\sup_{s \in [0, \theta]} \left| \frac{\mathrm{d}}{\mathrm{d}s} P_i(s) \right| < \infty.
\]
Hence the requirements of Definition \ref{def.contraction} of an RSCC with contraction
are met with $j=1$. 
By Theorem \ref{th.DF} it follows that the Markov chain $(s_{n, a})_{n \in \mathbb{N}}$ associated with this RSCC
with contraction is a Doeblin-Fortet chain and its transition operator $U$ is a Doeblin-Fortet operator.
It remains to prove the regularity of $U$ with respect to $L([0, \theta])$. 
For this we have to prove the existence of a point $s^{*} \in [0, \theta]$ such that 
$\lim_{n \rightarrow \infty} \left| \sigma_n(s) - s^{*} \right| = 0$, for any $s \in [0, \theta]$, 
where $\sigma_n(s)$ is the support of measure $Q^{n}(s, \cdot)$, $n \in \mathbb{N}_+$. 

Let $s \in [0, \theta]$ be an arbitrarily fixed number and define 
\begin{equation}
w_1 := s, \quad w_{n+1} := \frac{1}{w_n + m \theta}, \quad n \in \mathbb{N}_+. \label{5.11}
\end{equation}
We have $w_n \in [0, \theta]$ and letting $n \rightarrow \infty$ in (\ref{5.11})
we get 
\[
w_n \rightarrow s^{*} := \frac{-1+\sqrt{1+4\theta^2}}{2\theta}.
\]
Clearly, $w_{n+1} \in \sigma_1(w_n)$ and Lemma \ref{lem.supp} and an induction argument show that 
$w_{n} \in \sigma_n(s)$, $n \in \mathbb{N}_+$. Thus
\[
d\left(\sigma_n(s), s^{*}\right) \leq \left|w_n - s^{*}\right| \rightarrow 0, \quad n \rightarrow \infty
\]
where $d$ stands for the Euclidian distance on the line. 
Now, the regularity of $U$ with respect to $L([0,\theta])$ follows from Theorem \ref{th.5.9}.

From (\ref{5.13}) and Theorem \ref{th.ordered} there exist a stationary probability measure $Q^{\infty}$ 
and two constants $q<1$ and $K$ such that
\begin{equation}
\left\| U^{n}f - U^{\infty}f \right\|_L \leq K q^n \left\|f\right\|_L, \quad n \in \mathbb{N}_+, \ f \in L([0, \theta]) \label{5.12}
\end{equation}
where $U^{n}f$ is as in (\ref{5.20}) and
\begin{equation}
U^{\infty}f = \int_{0}^{\theta} f(x) Q^{\infty}(\mathrm{d}x). \label{5.200}
\end{equation}
Here $Q^{\infty} = \gamma_{\theta}$ is the invariant probability measure of the transformation $T_{\theta}$ in (\ref{1.7}),
i.e., $Q^{\infty}$ has the density $\rho_{\theta} = 1/(x + m \theta)$, $x \in [0, \theta]$, with the normalizing factor 
$1/\log(1+\theta^2)$.

\hfill $\Box$\\
\begin{remark}
Another way to put this is that $\rho_{\theta}$ is the eigenfunction of eigenvalue $1$ of the Perron-Frobenius operator $U$.
\end{remark}

\section{A Gauss-Kuzmin-type theorem}

Now, we may determine the limit of the sequence $(\mu(T^n_{\theta}<x))_{n \in {\mathbb{N}}_+}$ as $n \rightarrow \infty$ and give the rate of this convergence.

\begin{theorem} \rm{(A Gauss-Kuzmin-type theorem for $T_{\theta}$)} \label{G-K-L}
Let $([0,\theta],{\cal B}_{[0,\theta]}, T_{\theta})$ be as in Definition 1.1(i). 
\begin{enumerate}

\item[(i)]
For a probability measure $\mu$ on $([0,\theta],{\cal B}_{[0,\theta]})$,
let the assumption (A) as follows:
\[(A) \quad \mbox{ $\mu$ is non-atomic and has a Riemann-integrable density.}\] 
Then for any probability measure $\mu$ which satisfies (A),
the following holds:
 
\begin{equation}
\lim_{n \rightarrow \infty}\mu (T_{\theta}^n < x) 
= \frac{1}{\log (1+\theta^2)}\log ((m\theta+x)\theta), \quad x \in [0, \theta]. 
\label{6.1}
\end{equation}

\item[(ii)]
In addition to assumption of $\mu$ in (i), 
if the density of $[0, \theta]\ni x\mapsto \mu([0,x])$ is  Lipschitz continuous, 
then there exist two positive constants $q < 1$ and $K$ such that 
for any $x \in [0, \theta]$ and $n \in \mathbb{N}_+$, the following holds:
\begin{equation}
\lim_{n \rightarrow \infty}\mu (T_{\theta}^n < x)= 
\frac{1+\alpha q^n}{\log (1+\theta^2)}\log ((m\theta+x)\theta), 
\label{6.2}
\end{equation}
where $\alpha := \alpha(\mu, n, x)$ with $\left|\alpha\right| \leq K$.

As a consequence, the $n$-th error term $e_n(\theta,\mu;x)$ of the Gauss-Kuzmin problem
is obtained as follows:

\begin{equation}
e_{n}(\theta,\mu;x)=
\frac{\alpha q^n}{\log (1+\theta^2)}\log ((m\theta+x)\theta). \label{6.3}
\end{equation}
\end{enumerate}

\end{theorem}

\noindent \textbf{Proof.} 

Let $T_{\theta}$ be as in (\ref{1.7}).
By Proposition \ref{prop.4.1.}(ii), we have 
\begin{equation}
\mu \left((T_{\theta})^{-n}(A)\right) 
= \int_{A} U^nf_0(x) \rho_{\theta}(x)\mathrm{d}x \quad \mbox{for any } n \in \mathbb{N}, A \in {\cal B}_{[0,\theta]} \label{6.3'}
\end{equation}
where $f_0(x)= \frac{x\theta+1}{\theta}(d\mu / d\lambda_{\theta})(x)$ for $x \in [0, \theta]$. 
If $d\mu / d\lambda_{\theta} \in L([0, \theta])$, by (\ref{5.200}) we have
\begin{equation}
U^{\infty} f_0 = \int_{0}^{\theta} f_0(x)\,Q^{\infty}(\mathrm{d}x) = \int_{0}^{\theta} f_0(x)\,\gamma_{\theta}(\mathrm{d}x) =\frac{1}{\log(1+\theta^2)}. \label{6.4}
\end{equation} 

Taking into account (\ref{5.12}), there exist two constants $q<1$ and $K$ such that
\begin{equation}
\|U^n f_0 - U^{\infty} f_0\|_L \leq K q^n\left\|f_0\right\|_L, \quad n \in \mathbb{N}_+. \label{6.5}
\end{equation}
Furthermore, consider the Banach space $C([0, \theta])$ of all real-valued continuous functions on $[0, \theta]$ 
with the norm $\|f\| := \sup_{x \in [0, \theta]}|f(x)|$. Since $L([0, \theta])$ is a dense subspace of $C([0, \theta])$ we have
\begin{equation}
\lim_{n \rightarrow \infty} \|(U^n - U^{\infty})f\| 
= 0 \quad \mbox{for all } f \in C([0, \theta]). \label{6.6}
\end{equation}
Therefore, (\ref{6.6}) is valid for a measurable function $f_0$ which is $Q^{\infty}$-almost surely continuous, that is, for a Riemann-integrable function $f_0$. 
Thus, we have
\begin{eqnarray}
\lim_{n \rightarrow \infty} \mu \left(T^n_{\theta} < x\right) 
&=& \lim_{n \rightarrow \infty} \int_{0}^{x} U^nf_0(u) \rho_{\theta}(u) \,\mathrm{d}u  \\
&=& \frac{1}{\log(1+\theta^2)} \int_{0}^{x} \,\rho_{\theta}(u)\,\mathrm{d}u \\
&=& \frac{1}{\log(1+\theta^2)} \log ((m\theta+x)\theta). 
\end{eqnarray}
Hence (\ref{6.1}) is proved.
\hfill $\Box$
\\

\begin{remark}

Since the Lebesgue measure $\lambda$ satisfies assumptions in both (i) and (ii)
of Theorem \ref{G-K-L},
(\ref{6.1}) and (\ref{6.2}) hold for the case $\mu=\lambda$.
Hence Theorem \ref{G-K-L} gives
the solution of the Gauss-Kuzmin problem for the pair $(T_{\theta},\mu)$ 
instead of $(\tau,\lambda)$ in (\ref{1.3}).

\end{remark}

\begin{remark}
Until now, the estimate of the convergence rate remains an open question. 
To obtain a better estimate of the convergence rate involved, we may use a Wirsing-type approach as in \cite{Sebe-2005, Sebe-2010, W-1974}.
\end{remark}

\end{document}